\documentclass{gtmon_a}
\pdfoutput=1

\usepackage[all]{xy}
\usepackage{enumerate}

\proceedingstitle{Proceedings of the Nishida Fest (Kinosaki 2003)}
\conferencestart{28 July 2003}
\conferenceend{8 August 2003}
\conferencename{International Conference in Homotopy Theory}
\conferencelocation{Kinosaki, Japan}

\editor{Matthew Ando}
\givenname{Matthew}
\surname{Ando}

\editor{Norihiko Minami}
\givenname{Norihiko}
\surname{Minami}

\editor{Jack Morava}
\givenname{Jack}
\surname{Morava}

\editor{W Stephen Wilson}
\givenname{W Stephen}
\surname{Wilson}

\title{Goodwillie towers and chromatic homotopy: an overview}

\author{Nicholas J Kuhn}
\givenname{Nicholas J}
\surname{Kuhn}
\address{Department of Mathematics\\
University of Virginia\\\newline
Charlottesville VA 22903\\
USA}
\email{njk4x@virginia.edu}
\urladdr{}

\volumenumber{10}
\issuenumber{}
\publicationyear{2007}
\papernumber{15}
\startpage{245}
\endpage{279}

\doi{}
\MR{}
\Zbl{}

\arxivreference{math.AT/0410342}


\received{13 October 2004}
\revised{12 July 2005}
\accepted{}
\published{29 January 2007}
\publishedonline{29 January 2007}
\proposed{}
\seconded{}
\corresponding{}
\version{}

\makeatletter
\def\cnewtheorem#1[#2]#3{\newtheorem{#1}{#3}[section]
\expandafter\let\csname c@#1\endcsname\c@thm}

\let\xysavmatrix\xymatrix
\def\xymatrix{\disablesubscriptcorrection\xysavmatrix}
\AtBeginDocument{\let\bar\wbar\let\tilde\wtilde\let\hat\what}

\newcommand{\Id}{\mathrm{Id}}

\newtheorem{thm}{Theorem}[section]
\cnewtheorem{lem}[thm]{Lemma}
\cnewtheorem{cor}[thm]{Corollary}
\cnewtheorem{prop}[thm]{Proposition}
\cnewtheorem{conj}[thm]{Conjecture}

\theoremstyle{definition}

\cnewtheorem{defn}[thm]{Definition}
\cnewtheorem{defns}[thm]{Definitions}
\cnewtheorem{quest}[thm]{Question}
\cnewtheorem{notation}[thm]{Notation}
\cnewtheorem{exm}[thm]{Example}
\cnewtheorem{exs}[thm]{Examples}
\makeautorefname{exm}{Example}

\theoremstyle{remark}

\cnewtheorem{rem}[thm]{Remark}
\cnewtheorem{rems}[thm]{Remarks}

\makeatother  

\newcommand{\hocolim}{\operatorname*{hocolim}}
\newcommand{\holim}{\operatorname*{holim}}
\newcommand{\colim}{\operatorname*{colim}}
\newcommand{\hofib}{\operatorname*{hofib}}
\newcommand{\hocofib}{\operatorname*{hocofib}}

\newcommand{\Map}{\operatorname{Map}}
\newcommand{\MapS}{\operatorname{Map_{\mathcal S}}}
\newcommand{\MapT}{\operatorname{Map_{\mathcal T}}}

\newcommand{\A}{{\mathcal  A}}
\newcommand{\B}{{\mathcal  B}}
\renewcommand{\C}{{\mathcal  C}}
\renewcommand{\D}{{\mathcal  D}}
\newcommand{\Sp}{{\mathcal  S}}

\newcommand{\T}{{\mathcal  T}}

\newcommand{\I}{{\mathcal  I}}
\renewcommand{\J}{{\mathcal  J}}
\newcommand{\E}{{\mathcal  E}}

\newcommand{\X}{{\mathcal  X}}
\newcommand{\Pp}{{\mathcal  P}}

\newcommand{\Alg}{{\mathcal Alg}}

\newcommand{\Module}{{\mathcal Mod}}

\newcommand{\PP}{{\mathbb P}}

\newcommand{\Sinfty}{\Sigma^{\infty}}
\newcommand{\Oinfty}{\Omega^{\infty}}

\newcommand{\sm}{\wedge}

\newcommand{\ra}{\rightarrow}
\newcommand{\xra}{\xrightarrow}

\newcommand{\xla}{\xleftarrow}

\newcommand{\tr}{\mathit{tr}}

\begin{document}

\begin{abstract}
This paper is based on talks I gave in Nagoya and Kinosaki in August
of 2003.  I survey, from my own perspective, Goodwillie's work on towers
associated to continuous functors between topological model categories,
and then include a discussion of applications to periodic homotopy as
in my work and the work of Arone--Mahowald.
\end{abstract}

\maketitle

\section{Introduction} \label{introduction}

About two decades ago, Tom Goodwillie began formulating his calculus of
homotopy functors as a way to organize and understand arguments being
used by him and others in algebraic $K$--theory.  Though it was clear
early on that his general theory offered a new approach to the concerns
of classical homotopy, and often shed light on older approaches, it is
relatively recently that its promise has been begun to be realized. This
has been helped by the recent publication of the last of Goodwillie's
series \cite{goodwillie1,goodwillie2,goodwillie3}, and by the support of
many timely new results in homotopical algebra and localization theory
allowing his ideas to be applied more widely.

At the Workshop in Algebraic Topology held in Nagoya in August 2002,
I gave a series of three talks entitled `Goodwillie towers: key features
and examples', in which I reviewed the aspects of Goodwillie's work that
I find most compelling for homotopy theory.  A first goal of this paper
is to offer a written account of my talks.  As in my talks, I focus on
towers associated to functors, ie the material of \cite{goodwillie3}. As
`added value' in this written version, I include some fairly extensive
comments about the general model category requirements for running
Goodwillie's arguments.

At the Conference on Algebraic Topology held in Kinosaki just previous
to the workshop, I discussed a result of mine \cite{kuhntate} that
says that Goodwillie towers of functors of spectra split after periodic
localization.  This is one of a number of ways discovered so far in which
Goodwillie calculus interacts beautifully with  homotopy as organized
by the chromatic point of view; another is the theorem of 
Greg Arone and Mark Mahowald \cite{am}.  A second goal is to survey
these results as well, and point to directions for the future.

The paper is organized as follows.

In \fullref{section 2}, I describe the major properties of Goodwillie
towers associated to continuous functors from one topological model
category to another. In \fullref{model cat section}, I discuss model
category prerequisites.  The basic facts about cubical diagrams
and polynomial functors are reviewed in \fullref{cubes section}.
The construction of the Goodwillie tower of a functor is given in
\fullref{tower construction}, and I sketch the main ideas behind the
proofs that towers have the properties described in \fullref{section 2}.

In \fullref{examples}, I discuss some of my favorite examples: Arone's
model for the tower of the functor sending a space $X$ to $\Sinfty
\Map(K,Z)$ \cite{arone}, the tower for the functor sending a spectrum $X$
to $\Sinfty \Oinfty X$, the tower of the identity functor on the category
of commutative augmented $S$--algebras, and tower for the identity functor
on the category of topological spaces as analyzed by 
Johnson,
Arone, Mahowald and 
Dwyer
\cite{bjohnson,am,aronedwyer}. Besides
organizing these in a way that I hope readers will find helpful,
I have also included some remarks that haven't appeared elsewhere,
eg I note (in \fullref{kp example}) that the bottom of the tower
for $\Sinfty \Oinfty X$ can be used to prove the Kahn--Priddy Theorem,
`up to one loop'.

The long \fullref{periodic section} begins with a discussion of how
Goodwillie towers interact with Bousfield localization.  Included is a
simple example (see \fullref{LE example}) that shows that the composite
of homogeneous functors between spectra need not again be homogeneous.  In
the remainder of the section, I survey three striking results in which the
Goodwillie towers discussed in \fullref{examples} interact with chromatic
homotopy theory: my theorems on splitting localized towers \cite{kuhntate}
and calculating the Morava K--theories of infinite loopspaces \cite{k6},
and Arone and Mahowald's work on calculating the unstable $v_n$--periodic
homotopy groups of spheres \cite{am}. All of these relate to telescopic
functors $\Phi_n$ from spaces to spectra constructed a while ago by Pete
Bousfield and I \cite{bousfield3,k2,bousfield7} using the Nilpotence and
Periodicity Theorems of Devinatz, Hopkins and Smith
\cite{dhs,hs}. This suggests that Goodwillie calculus
can be used to further explore these curious functors.  Included in this
section, as an application of my work in \cite{k6}, is an outline of a
new way to possibly find a counterexample to the Telescope Conjecture.

The Kinosaki conference was on the occasion of Professor Nishida's
60th birthday, and I wish to both offer him my hearty congratulations,
and thank him for his kind interest in my research over the years.
Many thanks also to Noriko Minami and the other conference organizers
for their hospitality.

This research was partially supported by a grant from the National
Science Foundation.

\section{Properties of Goodwillie towers} \label{section 2}

The basic problem that Goodwillie calculus is designed to attack is
as follows.  One has a homotopy functor
$$F \co \C \ra \D$$
between two categories in which one can do homotopy. One wishes to
understand the homotopy type of $F(X)$, perhaps for some particular $X
\in \C$.

Goodwillie's key idea is to use the {\em functoriality} as $X$ varies,
to construct a canonical {\em polynomial} resolution of $F(X)$ as a
functor of $X$.

The first thing to specify is what is meant by categories in which one
can do homotopy theory.  In Goodwillie's papers, these are $\T$, the
category of pointed topological spaces, or $\Sp$, an associated category
of spectra (eg the $S$--modules of 
Elmendorf--Kriz--Mandell--May \cite{ekmm}), or variants of these,
eg $\T_Y$, the category of spaces over and under a fixed space $Y$.
But the arguments and constructions of Goodwillie
\cite{goodwillie3} are written
in a such a manner that they apply to situations in which $\C$ and $\D$
are suitably nice based model categories: in \fullref{model cat section},
we will spell out precisely what we mean.

Among all functors $F \co \C \ra \D$, some will be {\em $d$--excisive}
(or {\em polynomial of degree at most $d$}).  This will be carefully
explained in \fullref{polynomial functors}; we note that a $0$--excisive
functor is one that is homotopically constant, a functor is $1$--excisive
if it sends homotopy pushout squares to homotopy pullback squares,
and a $(d-1)$--excisive functor is also $d$--excisive.

Goodwillie's first theorem says that any $F$ admits a canonical polynomial
resolution.

\begin{thm} {\rm \cite[Theorem 1.8]{goodwillie3}}\qua \label{theorem 1}  Given
a homotopy functor $F \co \C \ra \D$ there exists a natural tower of
fibrations under $F(X)$,
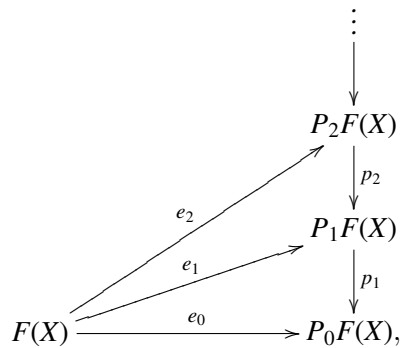
\begin{equation*}
\xymatrix{
&&& \vdots \ar[d] \\
&&& P_2F(X) \ar[d]^{p_2} \\
&&& P_1F(X) \ar[d]^{p_1} \\
F(X) \ar[rrr]^{e_0}  \ar[urrr]^{e_1} \ar[uurrr]^{e_2} &&& P_0F(X),
}
\end{equation*}
such that 
\begin{enumerate}
\item $P_dF$ is $d$--excisive, and \
\item $e_d \co F \ra P_dF$ is the universal weak natural
transformation to a $d$--excisive functor.
\end{enumerate}
\end{thm}

Let us explain what we mean by property (2).  By a weak natural
transformation $f \co F \ra G$, we mean a pair of natural transformations
$\smash{F\xra{g}H\xla{h}G}$ such that $\smash{H(X)\xla{h}G(X)}$ is a
weak equivalence
for all $X$.  Note that a weak natural transformation induces a well
defined natural transformation between  functors taking values in
the associated homotopy category.  Property (2) means that, given any
$d$--excisive functor $G$, and natural transformation $f \co F \ra G$,
there exists a weak natural transformation $g \co P_dF \ra G$ such that,
in the homotopy category of $\D$,
\begin{equation*}
\xymatrix{
F(X) \ar[d]_-{f(X)} \ar[r]^-{e_d(X)} & P_dF(X) \ar[dl]^-{g(X)}  \\
G(X)  & }
\end{equation*}
commutes for all $X \in \C$, and any two such $g$ agree.

A very useful property of the $P_d$ construction is the following Lemma.

\begin{lem} \label{fiber lem} Given natural transformations $F \ra G
\ra H$, if
$$ F(X) \ra G(X) \ra H(X)$$
is a fiber sequence for all $X$, then so is
$$ P_dF(X) \ra P_dG(X) \ra P_dH(X).$$
\end{lem}

Let $D_dF \co \C \ra \D$ be defined by letting $D_dF(X)$ be the homotopy
fiber of $P_dF(X) \ra P_{d-1}F(X)$.  The lemma and theorem formally imply
that $D_dF$ is {\em homogeneous of degree $d$}: it is $d$--excisive,
and $P_{d-1}D_dF(X) \simeq *$ for all $X$.

When $\D$ is $\T$, Goodwillie discovered a remarkable fact: these fibers
are canonically infinite loopspaces.  For a general $\D$, we let $\Sp(\D)$
be the associated category of `$\D$--spectra' (see \fullref{model cat
section}), and Goodwillie's second theorem then goes as follows.

\begin{thm} {\rm \cite[Theorem 2.1]{goodwillie3}}  \qua \label{theorem 2}  
Let $F \co \C \ra \D$ be homogeneous of degree $d$.  Then there is a naturally
defined homogeneous degree $d$ functor $F^{st} \co \C \ra \Sp(\D)$, such
that, for all $X \in \C$, there is a weak equivalence
$$ F(X) \simeq \Oinfty(F^{st}(X)).$$
\end{thm}

The category $\Sp(\D)$ is an example of a {\em stable} model category.
In a manner similar to results in the algebra literature, Goodwillie
relates homogenous degree $d$ functors landing in a stable model category
to symmetric multilinear ones.  A functor $L \co \C^d \ra \D$ is $d$--linear
if it is homogeneous of degree 1 in each variable, and is symmetric
if $L$ is invariant under permutations of the coordinates of $\C^d$.
Goodwillie's third theorem goes as follows.

\begin{thm} {\rm \cite[Theorem 3.5]{goodwillie3}} \qua \label{theorem 3}  Let 
$F \co \C \ra \D$ be a homogeneous functor of degree $d$ with $\D$ a stable model
category.  Then there is a naturally defined symmetric $d$--linear
functor $LF \co \C^d \ra \D$, and a weak natural equivalence
$$ (LF(X, \dots, X))_{h\Sigma_d} \simeq F(X).$$
\end{thm}

If $F \co \C \ra \D$ is a homotopy functor with $\C$ and $\D$ either $\T$
or $\Sp$, let $C_F(d) = L(D_dF)^{st}(S, \dots, S)$, a spectrum with
$\Sigma_d$--action. (Here $S$ denotes either the 0--sphere $S^0$ or its
suspension spectrum.)  Goodwillie refers to $C_F(d)$ as the $d^{th}$
Taylor coefficient of $F$ due the following corollary of the last theorem.

\begin{cor} \label{thm 3 cor}
In this situation, there is a weak natural transformation
$$ (C_F(d) \sm X^{\sm d})_{h\Sigma_d} \ra (D_dF)^{st}(X)  $$
that is an equivalence if either $X$ is a finite complex, or $F$ commutes
with directed homotopy colimits up to weak equivalence.
\end{cor}

As will be illustrated in the examples, these equivariant spectra have
often been identified.

The theorems above are the ones I wish to stress in these notes,
but I should say a little about convergence.  In \cite{goodwillie2},
Goodwillie carefully proves a generalized Blakers--Massey theorem, and
uses it to study questions that are equivalent to the convergence of
these towers in the cases when $\C$ is $\T$ or $\Sp$.  In particular,
many functors can be shown to be `analytic', and an analytic functor $F$
admits a `radius of convergence' $r(F)$ with the property that the tower
for $F(X)$ converges strongly for all $r(F)$--connected objects $X$.
The number $r(F)$ is often known, as will be illustrated in the examples.
A nice result from \cite{goodwillie2} reads as follows.

\begin{prop} {\rm \cite[Proposition 5.1]{goodwillie2}} \qua Let $F \ra G$ 
be a natural
transformation between analytic functors, and let $r$ be the maximum
of $r(F)$ and $r(G)$. If $F(X) \ra G(X)$ is an equivalence for all $X$
that are equivalent to high suspensions then it is an equivalence for
all $r$--connected $X$.
\end{prop}

\section{Model category prerequisites} \label{model cat section}

References for model categories include Quillen's original 1967
lecture notes \cite{quillen}, Dwyer and Spalinski's 1995 survey
article \cite{ds}, and the more recent books by Hovey and Hirschhorn
\cite{hovey1,hirschhorn}.

\subsection{Nice model categories} \label{nice model categories}

We will assume that $\C$ and $\D$ are either {\em simplicial} or {\em
topological} based model categories. (See \cite[Chapter 9]{hirschhorn}
for a careful discussion of simplicial model categories and \cite[Chapter
VII]{ekmm} for topological model categories.) `Based' means that the
initial and final object are the same: we will denote this object by $*$.

As part of the structure of a based topological (or simplicial) model
category $\C$, given $K \in \T$ and $X \in \C$, one has new objects
in $\C$, $X \otimes K$ and $\Map(K,X)$ satisfying standard properties.
This implies that $\C$ supports canonical homotopy limits and colimits:
given a functor $\X: \J \ra \C$ from a small category $\J$, $\hocolim_J
\X$ and $\holim_J \X$ are defined as appropriate coends and ends:
$$\hocolim_J \X = \X \otimes_{\J} E\J^{op}
  \quad\text{and}\quad \holim_J \X  = \Map^{\J}(E\J, \X),$$
where $E\J \co \J \ra \T$ is given by $E\J(j) = B(\J\downarrow j)_{+}$. (Thus
these generalize the constructions $X \times_G EG$ and $\Map^G(EG,X)$
familiar from equivariant topology.)  For details and basic properties,
see \cite[Chapter 18]{hirschhorn}.

With such canonical homotopy limits and colimits, $\C$ will support a
sensible theory of homotopy Cartesian and coCartesian cubes, as discussed
in \cite{goodwillie2}: see \fullref{cubes} below.  To know that certain
explicit cubes in $\C$ are homotopy coCartesian, one also needs that
$\C$ be left proper, and it seems prudent to require both $\C$ and $\D$
to be {\em proper}: the pushout of a weak equivalence by a cofibration
is a weak equivalence, and dually for pullbacks.

$\D$ then needs a further axiom ensuring that the sequential homotopy
colimit of homotopy Cartesian cubes is again homotopy Cartesian: assuming
that $\D$ admits the (sequential) {\em small object argument} does the
job: see Schwede \cite[Section 1.3]{schwede}.

\begin{exs}  The following categories satisfy our hypotheses:
\begin{itemize}
\item $\T_Y$, the category of spaces over and under $Y$,
\item $R-\Module$, the category of $R$--modules, where $R$ is an
$E_{\infty}$ ring spectrum, a.k.a. commutative $S$--algebra \cite{ekmm},
\item $R-\Alg$, the category of augmented commutative $R$--algebras,
\item simplicial versions of all of these, eg spectra as in \cite{bf}.
\end{itemize}
\end{exs}

\subsection{Spectra in model categories} \label{spectra}  Let $\D$ be
a model category as above, and let $\Sigma X$ denote $X \otimes S^1$,
ie $X$ tensored with the circle $S^1$.  Trying to force the suspension
$\Sigma \co \D \ra \D$ to be `homotopy invertible'  leads to a model category
of spectra $\Sp(\D)$ in the `usual way': this has been studied carefully
by Schwede \cite{schwede} (following Bousfield--Friedlander \cite{bf}), 
Hovey \cite{hovey2},
and Basterra--Mandell \cite{basterramandell}.  Roughly put, an object in
$\Sp(\D)$ will consist of a sequence of objects $X_0, X_1, X_2, \dots$
in $\D$, together with maps $\Sigma X_n \ra X_{n+1}$.  The point of
this construction is that the model category structure $\Sp(\D)$ has
the additional property that it is {\em stable}: homotopy cofibration
sequences in $\Sp(\D)$ agree with the homotopy fibration sequences.
The associated homotopy category will be triangulated.

As in the familiar case when $\D = \T$, there are adjoint functors
$$ \Sinfty \co \D \ra \Sp(\D) \text{ and } \Oinfty \co \Sp(\D) \ra \D.$$

If $\D$ is already stable these functors form a Quillen equivalence. For
an arbitrary $\D$, this adjoint pair can take a surprising form, as the
following example illustrates.

\begin{exm}  In \cite{basterramandell}, the authors show that the category
$\Sp(R{-}\Alg)$ is Quillen equivalent to $R-\Module$ so that $\Sinfty
\co R{-}\Alg \ra \Sp(R{-}\Alg)$ identifies with the Topological Andr\'e--Quillen
Homology functor\footnote{To be precise, by $TAQ(B)$ we mean the
Topological Andr\'e--Quillen Homology of $B$ with coefficients in the
$B$--bimodule $R$.}  $TAQ\co R{-}\Alg \ra R{-}\Module$, and $\Oinfty\co
\Sp(R{-}\Alg) \ra R{-}\Alg$ identifies with the functor sending an $R$--module
$M$ to the trivial augmented $R$--algebra $R \vee M$.  (Partial results
along these lines were also proved in \cite{basterramccarthy,schwede}.)
\end{exm}

\subsection{Functors between model categories} \label{functors}  Suppose
$\C$ and $\D$ are nice topological model categories.  There are couple
of useful properties that a functor
$$ F \co \C \ra \D$$
might have.

Firstly $F$ will usually be {\em continuous}:  for all $X$ and $Y$
in $\C$, the function
$$F \co \Map_{\C}(X,Y) \ra \Map_{\D}(F(X),F(Y))$$
should be continuous.

If $F$ is continuous, given $X \in \C$ and $K \in \T$, there is a natural
assembly map
\begin{equation} \label{assembly map}
F(X) \otimes K \ra F(X \otimes K)
\end{equation}
defined by means of various adjunctions.
The existence of these assembly maps implies that $F$ will be a {\em
homotopy} functor: a weak equivalence between fibrant cofibrant objects
in $\C$ is carried by $F$ to a weak equivalence in $\D$.

The second property that some functors $F$ satisfy is that $F$ commutes
with filtered homotopy colimits, up to weak equivalence.  A functor
having this property has sometimes also been termed `continuous', but
Goodwillie \cite{goodwillie3} more cautiously uses the term {\em finitary}
and so will we.

There are many interesting functors that are not finitary, as the next
example shows.

\begin{exm}  Let $L_E \co \Sp \ra \Sp$ be Bousfield localization of spectra
with respect to a spectrum $E$.  Then $L_E$ is finitary exactly when
the assembly map
$$ L_E(S) \sm X \ra L_E(X)$$
is a weak equivalence for all spectra $X$.  In other words, $L_E$ is
finitary exactly when it is smashing, a property that many interesting
$L_E$'s do not have.

Just to confuse the issue, we note that if $L_E$ is regarded as
taking values in the topological model category $L_E\Sp$, in which
equivalences are $E_*$--isomorphisms and fibrant objects are $E_*$--local
\cite[Chapter VIII]{ekmm}, then $L_E \co \Sp \ra L_E\Sp$ {\em is} finitary.
\end{exm}

Finally, lets say a word about maps between functors.  If $\C$ is not
small, then it seems a bit daunting (set theoretically) to impose a
model category structure on the class of functors $F \co \C \ra \D$.  As an
adequate fix for calculus purposes, we use the following terminology.
Call a natural transformation $f \co F \ra G$ a weak equivalence, and write
$F\smash{\xra[{\mbox{$\sim$}}]{f}} G$, if $f(X) \co F(X) \ra G(X)$ is a weak equivalence for
all $X$ in $\C$.  By a weak natural transformation $f \co F \ra G$ we mean a
pair of natural transformations of the form $F
\smash{\xla[\mbox{$\sim$}]{g}} H \smash{\xra{h}} G$
or $F \smash{\xra{h}} H \smash{\xla[\mbox{$\sim$}]{g}} G$.  We say that a diagram of weak natural
transformations commutes if, after evaluation on any object $X$, the
associated diagram commutes in the homotopy category of $\D$.  Finally,
we say that a diagram of functors $F \ra G \ra H$ a fiber sequence if
$F(X) \ra G(X) \ra H(X)$ is a (homotopy) fiber sequence for all $X$.

\section{Cubical diagrams and polynomial functors} \label{cubes section}

\subsection{Cubical diagrams} \label{cubes}

We review some of the theory of cubical diagrams; a reference is 
\cite[Section 1]{goodwillie2}.

Let $S$ be a finite set.  The power set of $S$, $\Pp(S) = \{T \subseteq
S\}$ , is a partially ordered set via inclusion, and is thus a small
category.  Let $\Pp_0(S) = \Pp(S) - \{\emptyset\}$ and let $\Pp_1(S) =
\Pp(S) - \{S\}$.

\begin{defns}
(a)\qua A {\em $d$--cube} in $\C$ is a functor 
$\X \co \Pp(S) \ra \C$ with $|S|=d$.

(b)\qua $\X$ is {\em Cartesian} if the natural map
$$ \X(\emptyset) \ra \holim_{T \in \Pp_0(S)} \X(T) $$
is a weak equivalence.

(c)\qua  $\X$ is {\em coCartesian} if the natural map
$$ \hocolim_{T \in \Pp_1(S)} \X(T) \ra \X(S) $$
is a weak equivalence.

(d)\qua  $\X$ is {\em strongly coCartesian} if $\X|_{\Pp(T)}:
\Pp(T) \ra \C$ is coCartesian for all $T \subseteq S$ with $|T| \geq 2$.

\end{defns}

Often $S$ will be the concrete set $\mathbf d = \{1,\dots,d\}$.

\begin{exm} A 0--cube $\X(0)$ is Cartesian if and only if it is coCartesian
if and only if $\X(0)$ is acyclic (ie weakly equivalent to the initial
object $*$).
\end{exm}

\begin{exm} \label{ex1} A 1--cube $f \co \X(0) \ra \X(1)$ is Cartesian if
and only if it is coCartesian if and only if $f$ is an equivalence.
\end{exm}

\begin{exm} \label{ex2} A 2--cube
\begin{equation*}
\xymatrix{
\X(0) \ar[d] \ar[r] & \X(1) \ar[d]  \\
\X(2) \ar[r] & \X(12) }
\end{equation*}
is Cartesian if it is a homotopy pullback square, and coCartesian if it
is a homotopy pushout square.
\end{exm}

\begin{exm} Strongly coCartesian $d$--cubes are equivalent to ones
constructed as follows.  Given a family of cofibrations $f(t) \co X(0) \ra
X(t)$ for $1 \leq t \leq d$, let $\X \co \mathbf d \ra \C$ be defined by
$\X(T) = \text{ the pushout of } \{f(t) \ | \ t\in T\}$.  (Note that
$\X(T)$ can be interpreted as the coproduct under $X(0)$ of $X(t),
t \in T$.)
\end{exm}

Critical to Goodwillie's constructions, is a special case of this last
example.

\begin{defn}  If $T$ is a finite set, and $X$ is an object in $\C$,
let $X*T$ be the homotopy cofiber of the folding map $\coprod_T X \ra X$.
\end{defn}

For $T \subseteq \mathbf d$, the assignment $T \mapsto X*T$ is easily
seen to define a strongly coCartesian $d$--cube $\X$: if $X \ra *$
factors as $\smash{X \xra{i} CX \xra{p} *}$, with $i$ a cofibration and $p$ an
acyclic fibration, then $\X$ agrees with the cube of the last example
with $f(t)=i \co X \ra CX$ for all $t$.

In the special case when $\C = \T$, $X*T$ is the (reduced) join of $X$
and $T$: the union of $|T|$ copies of the cone $CX$ glued together along
their common base $X$. 

There is a very useful way to inductively identify Cartesian cubes.
Note that the fibers of the vertical maps in a Cartesian 2--cube
as in \fullref{ex2} form a Cartesian 1--cube as in \fullref{ex1}.
This generalizes to higher dimensional cubes as we now explain.

Regard $\mathbf d$ as the obvious subset of $\mathbf {d+1}$.  Given an
$(d+1)$--cube $\X \co \Pp(\mathbf {d+1}) \ra \C$, we define three associated
$d$-cubes
$$\X_{\text{top}}, \X_{\text{bottom}}, \partial \X \co \Pp(\mathbf d)
\ra \C$$
as follows.  Let $ \X_{\text{top}}(T) = \X(T)$ and $\X_{\text{bottom}}(T)
= \X(T \cup \{n+1\})$.  Then define $\partial \X(T)$ by taking homotopy
fibers of the evident natural transformation between these:
$$ \partial \X(T) = \hofib \{\X_{\text{top}}(T) \ra
\X_{\text{bottom}}(T)\}.$$

\begin{lem} $\X$ is Cartesian if and only if $\partial \X$ is Cartesian.
\end{lem}

\begin{lem} If $\X_{\text{top}}$ and $\X_{\text{bottom}}$ are Cartesian,
so is $\X$.
\end{lem}

\begin{rem} Dual lemmas hold for coCartesian cubes.  One application of
this is that if $\C$ is a {\em stable} model category, so that homotopy
fibre sequences are the same as homotopy cofiber sequences, then $\X$
is Cartesian if and only if $\X$ is coCartesian.
\end{rem}

\subsection{Polynomial functors} \label{polynomial functors}

Let $\C$ and $\D$ be topological or simplicial model categories as in
\fullref{nice model categories}.

\begin{defn}  $F \co \C \ra \D$ is called {\em $d$--excisive} or said to be
{\em polynomial of degree at most $d$} if, whenever $\X$ is a strongly
coCartesian $(d+1)$--cube in $\C$, $F(\X)$ is a Cartesian cube in $\D$.
\end{defn}

\begin{exm}  $F$ has degree 0 if and only if $F(X) \ra F(*)$ is an
equivalence for all $X \in \C$, ie $F$ is homotopy constant.
\end{exm}

\begin{exm}  $F \co \C \ra \D$ is $1$--excisive means that $F$ takes pushout
squares to pullback squares.

In the classical case when $\C$ and $\D$ are spaces or spectra, this
implies that the functor sending $X$ to $\pi_*(F(X))$ satisfies the
Mayer--Vietoris property.

If $F$ is also finitary, then Milnor's wedge axiom holds as well. Then
there are spectra $C_0$ and $C_1$ such that $F(X) \simeq C_0 \vee (C_1
\sm X)$ if $\D = \Sp$ and $F(X) \simeq \Oinfty (C_0 \vee (C_1 \sm X))$
if $\D = \T$.
\end{exm}

\begin{rem} Without the finitary hypothesis, classifying $1$--excisive
functors seems very hard.  Examples of $1$--excisive functors of $X$
from spectra to spectra include the localization functors $L_EX$ and
functors of the form $\MapS(C,X)$ where $C \in \Sp$ is fixed.
\end{rem}

The following proposition of Goodwillie constructs $d$--excisive functors
out of $d$--variable $1$--excisive functors.

\begin{prop} {\rm \cite[Proposition 3.4]{goodwillie2}} \qua 
If $L \co \C^d \ra \D$ is
$1$--excisive in each of the $d$--variables, then the functor sending $X$
to $L(X, \dots, X)$ is $d$--excisive.
\end{prop}

\begin{cor} In this situation, if $L$ is symmetric, and $\D$ is a stable
model category, then, given any subgroup $G$ of the $d^{th}$ symmetric
group $\Sigma_n$, the  functor sending $X$ to $L(X, \dots, X)_{hG}$
is $d$--excisive.
\end{cor}

The various lemmas about identifying Cartesian cubes can be used to
prove the next two useful lemmas.

\begin{lem} If $F$ is $d$--excisive, then $F$ is $c$--excisive for all
$c \geq d$.
\end{lem}

\begin{lem} If $F \ra G \ra H$ is a fiber sequence of functors, and $G$
and $H$ are both $d$--excisive, then so is $F$.
\end{lem}

\section{Construction of Goodwillie towers and the proof of the \break main
properties} \label{tower construction}

\subsection[Construction of the tower and the proof of Theorem 1]
{Construction of the tower and the proof of \fullref{theorem 1}}

If one is to construct a $d$--excisive functor $P_dF$, then $P_dF(\X)$
needs to be Cartesian for for all strongly coCartesion $(d+1)$--cubes
$\X$. The idea behind the construction of $P_dF$ is to force this
condition to hold for certain strongly coCartesion $(d+1)$--cubes $\X$.

Fix an object $X \in \C$.  As discussed above, for $T \subseteq \mathbf
{d+1}$, the assignment $T \mapsto X*T$ defines a strongly coCartesian
$(d+1)$--cube $\X$.  For example, when $d+1 = 2$, one gets the pushout square
\begin{equation*}
\xymatrix{
X \ar[d] \ar[r] & CX \ar[d]  \\
CX \ar[r] & \Sigma X. }
\end{equation*}

\begin{defn} Let $T_dF \co \C \ra \D$ be defined by
$$T_dF(X) = \holim_{T \in \PP_0(\mathbf {d+1})} F(X*T).$$
\end{defn}

Note that there is an evident natural transformation $t_d(F) \co F \ra T_dF$,
and that this is an equivalence if $F$ is $d$--excisive.

\begin{defn} Let $P_dF \co \C \ra \D$ be defined by
$$ P_dF(X) = \hocolim \Bigl\{F(X) \xra{t_d(F)} T_dF(X) \xra{t_d(T_dF)}
T_dT_dF(X) \ra \cdots \Bigr\}.$$
\end{defn}

\begin{exm} $T_1F(X)$ is the homotopy pullback of
\begin{equation*}
\xymatrix{
 & F(CX) \ar[d]  \\
F(CX) \ar[r] & F(\Sigma X), }
\end{equation*}
Suppose that $F(*) \simeq *$.  Then $F(CX) \simeq *$, so that $T_1F(X)$
is equivalent to the homotopy pullback of
\begin{equation*}
\xymatrix{
 & {*} \ar[d]  \\
{*} \ar[r] & F(\Sigma X), }
\end{equation*}
which is $\Omega F(\Sigma X)$.  It follows that there is a natural
weak equivalence
$$P_1F(X) \simeq \hocolim_{n \ra \infty} \Omega^n F(\Sigma^n X).$$
\end{exm}

\begin{exm} \label{P1 example} Specializing the last example to the case
when $F$ is the identity functor $\Id \co \D \ra \D$, we see that
$$ P_1(\Id)(X) \simeq \Oinfty \Sinfty X.$$

If $\D = \T$, topological spaces, we see that $P_1(\Id)(X) = QX$.

If $\D = R-\Alg$, we see that $P_1(\Id)(B) \simeq R \vee TAQ(B)$ for an
augmented commutative $R$--algebra $B$.
\end{exm}

The proof of \fullref{theorem 1} amounts to checking that the $P_d$
construction just defined has the two desired properties: $P_dF$
should always be $d$--excisive, and $F \ra P_dF$ should be universal.
Checking the first of these is by far the more subtle, and follows from
the next lemma.

\begin{lem}[Goodwillie {{\cite[Lemma 1.9]{goodwillie3}}}] If $F \co \C \ra \D$ is a homotopy
functor, and $\X$ is strongly coCartesian $(d+1)$--cube in $\C$, then
there is a Cartesian $(d+1)$--cube $\mathcal Y$ in $\D$, such that $F(\X)
\ra T_dF(\X)$ factors through $\mathcal Y$.
\end{lem}

The construction of $\mathcal Y$ is very devious. $\mathcal Y$ is
(roughly) constructed to be the homotopy limit of $(d+1)$--cubes in $\D$
that are each seen to be Cartesian for the following reason: they are
constructed by applying $F$ to $(d+1)$--cubes in $\C$ formed by means
of evident objectwise equivalences between two $d$--cubes.

In contrast, proving that $F \ra P_dF$ is appropriately universal is
much easier.  Once one knows that $P_dF$ is $d$--excisive, universality
amounts to checking the following two things: 
\begin{enumerate}[(a)]
\item If $F$ is $d$--excisive, then $e_d(F) \co F \ra P_dF$ is a
weak equivalence, and 
\item $P_d(e_d(F))) \co P_dF \ra P_dP_dF$ is a weak equivalence. 
\end{enumerate}

These follow immediately from $T_d$--versions of these statements: 
\begin{enumerate}[(a${}^\prime$)]
\item If $F$ is $d$--excisive, then $t_d(F) \co F \ra T_dF$
is a weak equivalence, and 
\item $P_d(t_d(F))) \co P_dF \ra P_dT_dF$ is a weak equivalence.
\end{enumerate}

As was noted above, the first of these is clear.  The second admits
a fairly simple proof based on the commutativity of iterated homotopy
inverse limits.  Similar reasoning verifies the next lemma, which in
turn implies  \fullref{fiber lem}, which said that $P_d$ preserves
fiber sequences.

\begin{lem} Given natural transformations $F \ra G \ra H$, if
$$ F(X) \ra G(X) \ra H(X)$$
is a fiber sequence for all $X$, then so is
$$ T_dF(X) \ra T_dG(X) \ra T_dH(X).$$
\end{lem}

\subsection[Delooping homogeneous functors and Theorem 2]
{Delooping homogeneous functors and \fullref{theorem 2}}

The most surprising property of Goodwillie towers is stated in
\fullref{theorem 2}.  This says that, for $d>0$, homogeneous $d$--excisive
functors are infinitely deloopable.  To show this, Goodwillie proves
his beautiful key lemma, which says that $P_dF(X) \ra P_{d-1}F(X)$
is always a principal fibration if $F$ is {\em reduced}: $F(*) \simeq *$.

\begin{lem} \cite[Lemma 2.2]{goodwillie3} \label{key lemma} Let $d>0$,
and let $F \co \C \ra \D$ be a reduced functor.  There exists a homogeneous
degree $d$ functor $R_dF \co \C \ra \D$ fitting into a fiber sequence
of functors
$$ P_dF \ra P_{d-1}F \ra R_dF.$$
\end{lem}

Iteration of the $R_d$ construction leads to \fullref{theorem 2}: if $F$
is homogeneous of degree $d$, then we can let $F^{st}(X)$ be the spectrum
with $r^{th}$ space $R_d^rF(X)$.

The proof of \fullref{key lemma} is yet another clever manipulation of
categories related to cubes.  As an indication of how this might work,
we sketch how one can construct a homotopy pullback square
\begin{equation*}
\xymatrix{
T_dF(X) \ar[d] \ar[r] & K_dF(X) \ar[d]  \\
T_{d-1}F(X) \ar[r] & Q_dF(X) }
\end{equation*}
with $K_dF(X) \simeq *$, in the case when $d=2$.

One needs to look at how one passes from $\Pp_0(2)$ to $\Pp_0(3)$.
In pictures, $\Pp_0(2)$ looks like
\begin{equation*}
\xymatrix @-1.5pc {
& &&  1 \ar@{-}[dl]   \\
2 \ar@{-}[rr] & & 12, &
}
\end{equation*}
while $\Pp_0(3)$ looks like
\begin{equation*}
\xymatrix @-1.5pc {
& &&  1 \ar@{-}[dl] \ar@{-}[dd]  \\
2 \ar@{-}[dd] \ar@{-}[rr] & & 12 \ar@{-}[dd]&  \\
& 3 \ar@{-}[dl] \ar@{-}'[r][rr]  &&  13 \ar@{-}[dl]  \\
23 \ar@{-}[rr] & & 123.&
}
\end{equation*}

Now we decompose the poset $\Pp_0(3)$ as
\begin{equation} \label{decomposition}
\xymatrix{
\A \cap \B \ar[d] \ar[r] & \B \ar[d]  \\
\A \ar[r] & \A \cup \B = \Pp_0(3),}
\end{equation}
where $\A$ is
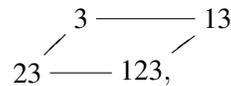
\begin{equation*}
\xymatrix @-1.5pc {
& 3 \ar@{-}[dl] \ar@{-}[rr]  &&  13 \ar@{-}[dl]  \\
23 \ar@{-}[rr] & & 123,&
}
\end{equation*}
and $\B$ is
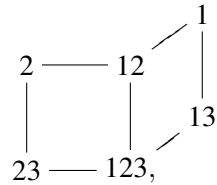
\begin{equation*}
\xymatrix @-1.5pc {
&  &&  1 \ar@{-}[dl] \ar@{-}[dd]  \\
2 \ar@{-}[dd] \ar@{-}[rr] & & 12 \ar@{-}[dd]&  \\
& &&  13 \ar@{-}[dl]  \\
23 \ar@{-}[rr] & & 123,&
}
\end{equation*}
so that $\A \cap \B$ is
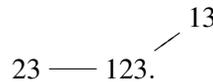
\begin{equation*}
\xymatrix @-1.5pc {
& &&  13 \ar@{-}[dl]  \\
23 \ar@{-}[rr] & & 123.&
}
\end{equation*}

The decomposition of posets \eqref{decomposition} induces a homotopy
pullback diagram
\begin{equation*}
\xymatrix{
\holim_{T \in \Pp_0(3)} F(X*T) \ar[d] \ar[r] & \holim_{T \in \A} F(X*T)
\ar[d]  \\
\holim_{T \in \B} F(X*T) \ar[r] & \holim_{T \in \A \cap \B} F(X*T). }
\end{equation*}

The top left corner is $T_2F(X)$, by definition.  As $\Pp_0(2)$ is cofinal
in $\B$, the bottom left corner is equivalent to $T_1F(X)$.  Finally $\A$
has initial object $\{3\}$, so that the upper left corner is contractible:
$$ \holim_{T \in \A} F(X*T) \simeq F(X * \{3\}) = F(CX) \simeq *.$$

\subsection[Cross effects and the proof of \ref{theorem 3}]
{Cross effects and the proof of \fullref{theorem 3}}

\begin{defn} Let $ F \co \C \ra \D$  be a functor. We define $cr_dF \co \C^d
\ra \D$, the $d^{th}$ {\em cross effect} of $F$, to be the the functor
of $d$ variables given by
\begin{equation*}
(cr_dF)(X_1, \dots, X_d) = \hofib \Bigl\{ F\Bigl(\bigvee_{i \in {\mathbf
d}} X_i\Bigr)
\ra \holim_{T \subset \Pp_0(d)} F\Bigl(\bigvee_{i \in {\mathbf d} - T}
X_i\Bigr)\Bigr\}.
\end{equation*}
\end{defn}

The $d$--cube sending $T$ to $\bigvee_{i \in {\mathbf d}
- T} X_i$ is easily seen to be strongly coCartesian; letting $d=2$
for example, the square
\begin{equation*}
\xymatrix{
X_1 \vee X_2 \ar[d] \ar[r] & X_2 \ar[d]  \\
X_1 \ar[r] & {*} }
\end{equation*}
is weakly equivalent to the evidently coCartesian square
\begin{equation*}
\xymatrix{
X_1 \vee X_2 \ar[d] \ar[r] & CX_1 \vee X_2 \ar[d]  \\
X_1 \vee CX_2 \ar[r] & CX_1 \vee CX_2. }
\end{equation*}
It follows that if $F$ is $(d-1)$--excisive, then $cr_dF(X_1, \dots, X_d)
\simeq *$ for all $X_i$.  A similar argument \cite[Lemma 3.3]{goodwillie3}
shows that if $F$ is $d$--excisive, then $cr_dF$ is $1$--excisive in
each of its variables.

More cleverly, Goodwillie establishes a converse of sorts when $\D$ is a
stable model category \cite[proof of Proposition 3.4]{goodwillie3}: if $F$
is $d$--excisive, and $cr_dF(X_1, \dots, X_d) \simeq *$ for all $X_i$,
then $F$ is $(d-1)$--excisive. This has the following useful consequence.

\begin{lem} \label{cross effects lemma} If $\D$ is a stable model
category, a natural map $\alpha: F \ra G$ between homogeneous functors
of degree $d$ will be an equivalence if and only if $cr_d\alpha \co cr_dF
\ra cr_dG$ is an equivalence.
\end{lem}

Another property of $cr_dF \co \C^d \ra \D$ that is easy to see is that it
is {\em reduced}: $$cr_dF(X_1, \dots, X_d) \simeq *$$ if $X_i \simeq *$
for some $i$.

A permutation of $\mathbf d$, $\sigma \in \Sigma_d$, induces an evident
isomorphism
$$ \sigma_* \co cr_dF(X_1, \dots, X_d) \ra cr_dF(X_{\sigma(1)}, \dots,
X_{\sigma(d)}), $$
satisfying $(\sigma \circ \tau)_* = \sigma_* \circ \tau_*$: a functor
of $d$ variables with this structure is called {\em symmetric}.

\begin{defn} Let $L_dF \co \C^d \ra \D$ be the functor obtained from $cr_dF$
by applying $P_1$ to each variable.  Thus we have
$$ L_dF(X_1, \dots, X_d) \simeq \hocolim_{n_i \ra \infty} \Omega^{n_1 +
\dots + n_d} cr_dF(\Sigma^{n_1}X_1, \dots, \Sigma^{n_d}X_d).$$
\end{defn}

$L_dF$ will always be symmetric and $d$--linear, and if $F$ is
$d$--excisive, then the natural map $cr_dF \ra L_dF$ is an equivalence.

If $G$ is a finite group, let $G-\D$ denote the category of objects
in $\D$ with a $G$--action.  Given $Y \in G-\D$, we let $Y_{hG} = Y
\otimes_G EG$ and $Y^{hG}= \Map^G(EG,Y)$ denote the associated homotopy
quotient and fixed point objects in $\D$.

\begin{defn} Let $\Delta_dF \co \C \ra \Sigma_d-\D$ be defined by
$$\Delta_dF(X) = L_dF(X, \dots, X).$$
\end{defn}

A more precise version of \fullref{theorem 3} is the following.

\begin{thm} \label{fiber formula theorem} Let $F \co \C \ra \D$ be a
homotopy functor, with $\D$ a stable model category.  Then there is a
natural weak equivalence
$$ \Delta_dF(X)_{h\Sigma_d} \simeq D_dF(X).$$
\end{thm}

If $F$ is $d$--excisive then $\Delta_dF(X)$ can be identified with
$(cr_dF)(X, \dots, X)$. In this case, one gets a natural transformation
$$ \alpha_d(X) \co (\Delta_dF)(X)_{h\Sigma_d} \ra F(X)$$
defined to be the composite
$$ (\Delta_dF)(X)_{h\Sigma_d} \ra F\Bigl(\bigvee_{i=1}^d X\Bigr)_{h\Sigma_d}
\ra F(X).$$
Here the second map is induced by the fold map $ \bigvee_{i=1}^d X
\ra X$.  The theorem is then proved by verifying that $cr_d(\alpha_d)$
is an equivalence, and appealing to \fullref{cross effects lemma}.

We indicate how \fullref{thm 3 cor} follows from \fullref{fiber formula
theorem}.  The assembly map for $F$ induces an assembly map
$$ (\Delta_dF(X) \otimes K^{\sm d})_{h \Sigma_d} \ra \Delta_dF(X \otimes
K)_{h \Sigma_d},$$
for $X \in \C$ and $K \in \T$.

If $\C$ is $\T$, $\D$ is $\Sp$, and $X = S$, then this reads
$$ (C_F(d) \sm K^{\sm d})_{h\Sigma_d} \ra \Delta_dF(K)_{h\Sigma_d},$$
where $C_F(d) = \Delta_dF(S)$.  By construction, this map is the identity
if $K = S$, and it follows that it will be an equivalence for all finite
$K$, or all $K$ under the additional hypothesis that $F$ is finitary.
A similar argument holds if both $\C$ and $\D$ are $\Sp$: here the
assembly map can be constructed for all $K \in \Sp$.

When the domain category $\C$ is also stable, and $\D = \Sp$, there is
an elegant addendum to \fullref{fiber formula theorem} essentially due
to R McCarthy \cite{mccarthy}.

Given $Y \in G-\Sp$, there is a natural norm map $N(Y) \co Y_{hG} \ra Y^{hG}$
satisfying the property that $N(Y)$ is an equivalence if $Y$ is a finite
free $G$--CW spectrum. As in \cite{kuhntate}, we let the Tate spectrum
of $Y$, $\T_G(Y)$, be the cofiber.

\begin{prop} {\rm \cite{kuhntate}}\qua \label{pullback prop}  Let 
$F \co \C \ra \Sp$
be any homotopy functor, with $\C$ stable.  For all $d \geq 1$, there
is a homotopy pullback diagram
\begin{equation*}
\xymatrix{
P_dF(X) \ar[d] \ar[r] &
(\Delta_dF(X))^{h\Sigma_d} \ar[d]  \\
P_{d-1}F(X) \ar[r] &
{\mathcal T}_{\Sigma_d}(\Delta_dF(X)).
}
\end{equation*}
\end{prop}

\section{Examples} \label{examples}

\subsection{Suspension spectra of mapping spaces}

Fix a finite C.W.~ complex $K$.  Let $\MapT(K,X)$ be the space of based
continuous maps from $K$ to a space $X$.  Similarly, given a spectrum $Y$,
let $\MapS(K,Y)$ be the evident function spectrum.

In \cite{goodwillie2}, Goodwillie proved that the functor from spaces
to spectra sending $X$ to $\Sinfty \MapT(K,X)$ is analytic with radius
of convergence equal to the dimension of $K$.

In \cite{arone}, Arone gave a very concrete model for the associated
Goodwillie tower $\{ P_*^K(X) \}$.   The paper by Ahearn and Kuhn
\cite{ahearnkuhn}
includes further details about Arone's construction while building in
extra structure.

Let $\E$ be the category with objects the finite sets $\mathbf d$, $d
\geq 1$, and with morphisms the epic functions.  $\E_d$ will denote the
full subcategory with objects $\mathbf c$ with $c \leq d$.

Given a based space $X$, let $X^{\sm} \co \E^{op} \ra \T$ be the functor
sending $\mathbf d$ to $X^{\sm d}$.  Then Arone's model for $P_d^K \co
\T \ra \Sp$ is given by
$$ P_d^K(X) = \MapS^{\E_d}(K^{\sm}, \Sinfty X^{\sm}),$$
the spectrum of natural transformations between the two contravariant
functors of $\E_d$.  The natural transformation
$$ \Sinfty \MapT(K,X) \ra P_d^K(X)$$
is induced by sending $f \co K \ra X$ to $f^{\sm} \co K^{\sm} \ra X^{\sm}$
and then stablizing.

A by product of this construction is that there is a homotopy
pullback square of $\Sp$--modules which has some of the same flavor as
\fullref{pullback prop}:
\begin{equation*} \label{pullback 1}
\xymatrix{
P_d^K(X) \ar[d] \ar[r] & \MapS^{\Sigma_d}(K^{\sm d}, \Sinfty X^{\sm d})
\ar[d]  \\
P_{d-1}^K(X) \ar[r] & \MapS^{\Sigma_d}(\delta_d(K), \Sinfty X^{\sm d}), }
\end{equation*}
where $\delta_d(K) \subset K^{\sm d}$ denotes the fat diagonal.

Thus the $\!d$th fiber, $\!D_d^K(X)$, can be described as follows.
Let $\!K^{(d)}$ denote $\!K^{\sm d}/\delta_d(K)$.  Then we have
\begin{equation*}
\begin{split}
D_d^K(X) & = \MapS^{\Sigma_d}(K^{(d)}, \Sinfty X^{\sm d}) \\
  & \simeq \MapS(K^{(d)}, \Sinfty X^{\sm d})_{h\Sigma_d} \\
  & \simeq (\D(K^{(d)}) \sm X^{\sm d})_{h\Sigma_d}.
\end{split}
\end{equation*}
Here $\D(K^{(d)})$ denotes the equivariant $S$--dual of $K^{(d)}$, and
the equivalences follow from the fact that $K^{(d)}$ is both finite and
$\Sigma_d$--free away from the basepoint. It follows that the $d$th
Taylor coefficient of the functor sending $X$ to $\Sinfty \MapT(K,X)$
is $\D(K^{(d)})$.

\begin{rem}  In \cite{k4}, we observed that, when $X$ is also a finite
complex, the tower $P^K_*(X)$ also arises as by taking the $S$--dual of
a natural filtration on the nonunital commutative $S$--algebra 
$\D(X) \otimes K$.
\end{rem}

By Alexander duality, $\D(K^{(d)})$ can be identified an appropriate
equivariant desuspension of the suspension spectrum of a configuration
space.
Specializing to the case when $K = S^n$, this takes the following concrete
form.  Let $\C(n,d)$ denote the space of $d$ distinct little $n$--cubes
in a big $n$--cube (May \cite{may}).  Via a Thom--Pontryagin collapse, there
is a very explicit duality map of $\Sigma_d$ spaces \cite{ahearnkuhn}
$$ \C(n,d)_+ \sm S^{n(d)} \ra S^{nd}.$$
One proof that Arone's model works when $K=S^n$ goes roughly as follows.
Suppose $X = \Sigma^n Y$.  One has the usual filtered configuration space
model $C_n(Y)$ for $\Omega^n \Sigma^n Y$ \cite{may}.  Thus one has maps
$$ \Sinfty F_dC_n(Y) \ra \Sinfty \Omega^n \Sigma^n Y \ra P^{S^n}_d(Y).$$
The nontriviality of the second map is proved by showing that the
composite is an equivalence.  By induction on $d$, it suffices to
show that $cr_d$ applied to this composite is an equivalence, and the
verification of that leads back to the above explicit duality map.

A bonus corollary of this proof is that one also establishes a rather
nice  version of `Snaith splitting': the tower strongly splits when $X =
\Sigma^n Y$.

\begin{exm}  One application of this comes from applying mod $p$
cohomology to the tower.  One obtains a spectral sequence of differential
graded algebras $\{E_r^{s,t}(S^n,X)\}$ with
$$ E_1^{d,*}(S^n,X) =  H^*((\C(n,d)_+ \sm (\Sigma^{-d}X)^{\sm d})_{h
\Sigma_d}; \mathbb Z/p))$$
and converging strongly to $H^*(\Omega^n X; \mathbb Z/p)$ if $X$ is
$n$--connected.  This $E_1$ term is a known functor of $H^*(X;\mathbb
Z/p)$. The differentials have not been fully explored, but seem to
be partly determined by derived functors of destablization of unstable
modules over the Steenrod algebra, as applied to the $\mathcal A$--module
$\Sigma^{-n}H^*(X; \mathbb Z/p)$.
\end{exm}

\subsection{Suspension spectra of infinite loopspaces}

The previous example can be used to determine the tower $\{
P_d^{S^{\infty}} \}$ for the functor from spectra to spectra sending
a spectrum $X$ to $\Sinfty \Oinfty X$.

Let $X_n$ denote the $n$th space of the spectrum $X$.  Then we have
that $ \Omega^n X_n \simeq \Oinfty X$ for all $n$, and the natural map
$$ \hocolim_{n \ra \infty} \Sigma^{-n}\Sinfty X_n \ra X$$
is an equivalence.  From this and the last example, one can deduce that
the tower converges for $0$--connected spectra $X$ and that
$$ \hocolim_{n \ra \infty} \Sigma^{-n} P_d^{S^{n}}(X_n) \simeq
P_d^{S^{\infty}}(X).$$

As $\hocolim_{n \ra \infty}  \C(n,d)_+$ is a model for
$E\Sigma_d+$, and this is weakly equivalent to $S^0$, it follows that
the formula for the $d$th fiber is
$$ D_d^{S^{\infty}}(X) \simeq X^{\sm d}_{h \Sigma_d},$$
and thus the $d$th Taylor coefficient of the functor sending a spectrum
$X$ to $\Sinfty \Oinfty X$ is the sphere spectrum $S$ for all $d >0$.

Finally, \fullref{pullback prop} specializes to say that for each $d>0$
there is a pullback square
\begin{equation*} \label{pullback 2}
\xymatrix{
P_d^{S^{\infty}}(X) \ar[d] \ar[r] & (X^{\sm d})^{h\Sigma_d} \ar[d]  \\
P_{d-1}^{S^{\infty}}(X) \ar[r] & \T_{\Sigma_d}(X^{\sm d}). }
\end{equation*}

\begin{exm} \label{kp example} The tower begins
\begin{equation*}
\xymatrix{
 & P_2^{S^{\infty}}(X) \ar[d]^{p_1}  \\
\Sinfty \Oinfty X \ar[ur]^{e_2} \ar[r]^{e_1} & X, }
\end{equation*}
where $e_1 \co \Sinfty \Oinfty X \ra X$ is adjoint to the identity on
$\Oinfty X$.
A formal consequence is that
$$\Oinfty p_1 \co \Oinfty P_2^{S^{\infty}}(X) \ra \Oinfty X$$
admits a natural section.

The map $p_1$ fits into a natural cofibration sequence
$$ P_2^{S^{\infty}}(X) \xra{p_1} X \ra \Sigma (X \sm X)_{h\Sigma_2}.$$
Specializing to the case when $X = S^{-1}$, this can be identified
\cite[Appendix]{kuhntate} with the cofibration sequence
$$ \Sigma^{-1}\mathbb RP^{\infty}_0 \xra{t} S^{-1} \ra \mathbb
RP^{\infty}_{-1},$$
where $t$ is one desuspension of the Kahn--Priddy transfer, and $\mathbb
RP^{\infty}_{k}$ denotes the Thom spectrum of $k$ copies of the canonical
line bundle over $\mathbb RP^{\infty}$.

Letting $QZ$ denote $\Oinfty \Sinfty Z$, we conclude that
$$ \Oinfty \tr \co \Omega Q\mathbb RP^{\infty}_+ \ra \Omega QS^0$$
admits a section: a result `one loop' away from the full strength of the
Kahn--Priddy Theorem \cite{kahnpriddy} at the prime 2.  The odd prime
version admits a similar proof using that, localized at a prime $p$,
$P_d^{S^{\infty}}(S^{-1}) \simeq *$ for $1<d<p$.
\end{exm}

\subsection[The identity functor for Alg]{ The identity functor for $\Alg$}

Let $\Alg$ be the category of commutative augmented $S$--algebras.
This is a model category in which weak equivalences and fibrations are
determined by forgetting down to $S$--modules.  More curious is that
the coproduct of $A$ and $B$ is $A \sm B$.

Let $\bigl\{\smash{P_d^{alg}}\bigr\}$ denote the tower associated to
the identity $\I \co \Alg \ra \Alg$.  Given $B \in \Alg$, it is not
too hard to deduce that the tower $\bigl\{\smash{P_d^{alg}}(B)\bigr\}$
will strongly converge to $B$ if $I(B)$ is 0--connected, where $I(B)$
denotes the `augmentation ideal': the homotopy fiber of the augmentation
$B \ra S$.

Let $\smash{D_d^{alg}}(B)$ be the fiber of $\smash{P_d^{alg}}(B)
\ra \smash{P_{d-1}^{alg}}(B)$.  As already discussed in \fullref{P1
example}, $\smash{D_1^{alg}}(B)$ can be identified with $TAQ(B)$, the
Topological Andr\'e--Quillen Homology of $B$ with coefficients in the
$B$--bimodule $S$.

The fact that coproducts in $\Alg$ correspond to smash products of
$S$--modules leads to a simple calculation of the $d^{th}$ cross effect
of $\I$:
$$ cr_d(\I)(B_1, \dots, B_d) \simeq I(B_1) \sm \cdots \sm I(B_d).$$
From this, one gets a formula for $D_d^{alg}(B)$:
\begin{thm}
$$ D_d^{alg}(B) \simeq TAQ(B)^{\sm d}_{h\Sigma_d}.$$
\end{thm}
A proof of this in the spirit of this paper appears in \cite{k6}.
See also Minasian \cite{minasian}.

A nice corollary of this formula says the following.
\begin{cor} \label{taq cor}
If $A$ and $B$ in $\Alg$ have 0--connected augmentation ideals, then an
algebra map $f \co A \ra B$ is an equivalence if $TAQ(f)$ is.
\end{cor}

The converse of this corollary --- that $TAQ(f)$ is an equivalence if $f$
is --- is true even without connectivity hypotheses: see, eg \cite{k4}.
Without any hypotheses implying convergence, one has that if $TAQ(f)$ is
an equivalence, so is $\hat f \co \hat A \rightarrow \hat B$, where $\hat A$
denotes the homotopy inverse limit of the tower for $A$ \cite{k6}.

\begin{exm}  \label{splitting example} This tower overlaps in an
interesting way with the one for $\Sinfty \Oinfty X$ discussed above, and
the corollary leads to a simple proof of a highly structured version of the
classical stable splitting \cite{kahn} of $QZ$, for a connected space $Z$.

It is well known that $\Sinfty (\Oinfty X)_+$ is an $E_{\infty}$ ring
spectrum.  Otherwise put, we can regard $\Sinfty (\Oinfty X)_+$ as an
object in $\Alg$.  It is not hard to see that $TAQ(\Sinfty (\Oinfty X)_+)$
is equivalent to the connective cover of $X$, and there is an equivalence
$$ S \vee P^{S^{\infty}}_d(X) \simeq P_d^{alg}(\Sinfty (\Oinfty X)_+)$$
for connective spectra $X$.

Another object in $\Alg$ is $\mathbb P(X)$, the free commutative
$S$--algebra generated by $X$.  As an $S$--module,
$$\mathbb P(X) \simeq \bigvee_{d=0}^{\infty} X^{\sm d}_{h\Sigma_d},$$
and it is not hard to compute that $TAQ(\mathbb P(X)) \simeq X$.

The stable splitting of $Q Z$ gets proved as follows.  The inclusion
$$ \eta(Z) \co Z \ra Q Z$$
induces a natural map in $\Alg$
$$ s(Z) \co \mathbb P(\Sinfty Z) \ra \Sinfty (QZ)_+.$$
The construction of $s$ makes it quite easy to verify that $TAQ(s(Z)) \co
\Sinfty Z \ra \Sinfty Z$ is the identity.  The above corollary then
implies that $s(Z)$ is an equivalence in $\Alg$, and thus in $\Sp$,
for connected spaces $Z$.

A more detailed discussion of this appears in \cite{k6}.
\end{exm}

\subsection[The identity functor for T]{The identity functor for $\T$}

Let $\{P_d\}$ denote the tower of the identity functor on $\T$.  Given a
space $Z$, let $D_d(Z)$ denote the fiber of $P_d(Z) \ra P_{d-1}(Z)$,
and then let $D_d^{st}(Z) \in \Sp$ be the infinite delooping provided
by \fullref{theorem 2}.

Goodwillie's estimates \cite{goodwillie2} show that the tower $\{P_d(Z)\}$
will strongly converge to $Z$ when $Z$ is simply connected. With a
little more care, one can show that this is still true if $Z$ is just
nilpotent.\footnote{Thanks to Greg Arone for this observation.}  Thus,
for nilpotent $Z$, one gets a strongly convergent 2nd quadrant spectral
sequence converging to $\pi_*(Z)$, with $E^1_{-d,*} = \pi_*(D_d^{st}(Z))$.

Johnson \cite{bjohnson} and Arone with collaborators Mahowald and Dwyer
\cite{am,aronedwyer} have identified the spectra $D_d^{st}(Z)$ thusly:
$$ D_d^{st}(Z) \simeq (\D(\Sigma K_d) \sm Z^{\sm d})_{h\Sigma_d}, $$
where $K_d$ is the unreduced suspension of the classifying space of
the poset of nontrivial partitions of $\mathbf d$, and $\D$ denotes the
equivariant $S$--dual, as before.

Using this model, the discovery of Arone and Mahowald \cite{am} is that
when $Z$ is an odd dimensional sphere, these spectra are very special
spectra that were known previously.  To state the theorem, we need
some notation.

Let $p$ be a prime.  Let $m \rho_k$ denote the direct sum of $m$ copies
of the reduced real regular representation of $V_k = (\mathbb Z/p)^k$.
Then $GL_k(\mathbb Z/p)$ acts on the Thom space $(BV_k)^{m\rho_k}$.
Let $e_k \in \mathbb Z_{(p)}[GL_k(\mathbb Z/p)]$ be any idempotent in
the group ring representing the Steinberg module, and then let $L(k,m)$
be the associated stable summand of  $(BV_k)^{m\rho_k}$:
$$ L(k,m) = e_k (BV_k)^{m\rho_k}.$$
The spectra $L(k,0)$ and $L(k,1)$ agree with spectra called $M(k)$
and $L(k)$ in the literature: see eg \cite{mitchellpriddy,k1,kuhnpriddy}.

Collecting results from \cite{am} and 
\cite[Theorem 1.9, Corollary 9.6]{aronedwyer},
one has the following theorem.

\begin{thm} \label{Arone theorem}  Let $m$ be an odd natural number. 
\begin{enumerate}
\item $D_d^{st}(S^{m}) \simeq *$ if $d$ is not a power of a
prime.
\item Let $p$ be a prime. $D_{p^k}^{st}(S^{m}) \simeq
\Sigma^{m-k}L(k,m)$, and thus has (only) $p$--torsion homotopy if $k>0$.
\item $H^*(L(k,m); \mathbb Z/p)$ is free over the subalgebra
$\A(k-1)$ of the Steenrod algebra.  As a function of $k$, the connectivity
of $L(k,m)$ has a growth rate like $p^k$.
\end{enumerate}
\end{thm}

Thus the associated spectral sequences for computing the unstable homotopy
groups of odd spheres coverges exponentially quickly, and begins from
stable information about spectra of roughly the same complexity as
the suspension spectra of classifying spaces of elementary abelian
$p$--groups.

\begin{rem}  When $m = 1$, one gets a spectral sequence converging to
the known graded group $\pi_*(S^1)$, with $E^1_{-k,*} = \pi_*(L(k))$.
Comparison with my work on the Whitehead Conjecture \cite{k1,kuhnpriddy}
suggests that $E^2 = E^{\infty}$.  Greg Arone and I certainly believe
this, but a rigorous proof has yet to be nailed down.
\end{rem}

As discussed near the end of the next section, the properties listed in
the theorem have particularly beautiful consequences for computing the
periodic unstable homotopy groups of odd dimensional spheres.

\section{Interactions with periodic homotopy} \label{periodic section}

For topologists who study classical unstable and stable homotopy theory,
a major development of the past two decades has been the organization
of these subjects via the chromatic filtration associated to the Morava
K--theories.

One of the most unexpected aspects of Goodwillie towers is that they
interact with the chromatic aspects of homotopy in striking ways.
In this section, I survey, in inverse order of when they were proved,
three different theorems of this sort.

\subsection{Goodwillie towers and homology isomorphisms}  There are a
couple of useful general facts about how Bousfield localization relates
to Goodwillie towers.

Let $E_*$ be a generalized homology theory.  A map $f \co X \ra Y$ of
spaces or spectra is called an {\em $E_*$--isomorphism} if $E_*(f)$ is
an isomorphism.  A natural transformation $f \co F \ra G$ between functors
$F,G \co \C \ra \Sp$ is an {\em $E_*$--isomorphism} if $f(X)$ is for all
$X \in \C$.  Then we have

\begin{prop}{\rm\cite[Corollary 2.4]{k2004}}\qua
\label{homology iso prop1} If $F:
\C \ra \Sp$ is finitary and $f \co X \ra Y$ is an $E_*$--isomorphism then
so are $D_dF(f) \co D_dF(X) \ra D_dF(Y)$ and $P_dF(f) \co P_dF(X) \ra P_dF(Y)$
for all $d$.
\end{prop}

\begin{prop}[{{\cite[Lemma 6.1]{kuhntate}}}]
\label{homology iso prop2}
If a natural transformation $f \co F \ra G$ between functors $F,G \co \C
\ra \Sp$ is an $E_*$--isomorphism then so are $D_df \co D_dF \ra D_dG$
and $P_df \co P_dF \ra P_dG$ for all $d$.
\end{prop}

Both of these follow by observing that the various constructions defining
$P_d$ and $D_d$ preserve $E_*$--isomorphisms.

The next example illustrates that the finitary hypothesis in
\fullref{homology iso prop1} is needed.

\begin{exm} Consider $L_{H\mathbb Z/p} \co \Sp \ra \Sp$, which is
a homogeneous
functor of degree 1.  Then $H\mathbb Z/p^{\infty}$ is $H\mathbb
Q_*$--acyclic (ie $H\mathbb Z/p^{\infty} \ra *$ is an $H\mathbb
Q_*$--equivalence), but $L_{H\mathbb Z/p}(H\mathbb Z/p^{\infty}) =
\Sigma H\mathbb Z_p$ is not.
\end{exm}

For an application of \fullref{homology iso prop1} to the homology
of mapping spaces, see \cite{k2004}.  \fullref{homology iso prop2}
is crucially used in proofs of the two theorems discussed in the next
two subsections.

We end this subsection with some observations related to the phenomenon
illustrated in the last example.

If $F \co \C \ra \Sp$ is homogeneous of degree $d$, the functor $L_EF \co
\C \ra \Sp$ will always again be $d$--excisive, but need no longer
be homogeneous.

\begin{exm} \label{LE example} Let $F \co \Sp \ra \Sp$ be defined by
$F(X) = (X \sm X)_{h\Sigma_2}$.  The composite functor $L_EF$ will be
$2$--excisive, but need no be longer homogenous, even when restricted
to finite spectra.  Indeed, a simple calculation shows that
$$ P_1(L_EF)(S) = \hocofib \{L_ES \sm \mathbb RP^{\infty} \ra L_E\mathbb
RP^{\infty} \}.$$
This can easily be nonzero.  For example, when $E$ is mod 2 $K$--theory,
one has that $L_E\mathbb RP^{\infty} = L_ES$, as the transfer $\mathbb
RP^{\infty} \ra S$ is an $K\mathbb Z/2_*$--isomorphism.  It follows that
$P_1(L_EF)(S)$ has nonzero rational homology.
\end{exm}

As a fix for this problem, we have the next lemma, which follows from
\fullref{homology iso prop2}.

\begin{lem} \label{homology error lemma} If $F \co \C \ra \Sp$ is homogeneous
of degree $d$, then $P_{d-1}(L_EF)$ is $E_*$--acyclic. Otherwise said,
$D_d(L_EF) \ra L_EF$ is an $E_*$--isomorphism.
\end{lem}

\subsection{Goodwillie towers and periodic localization}

We will consider two families of periodic homology theories.

Fixing a prime $p$,  $K(n)_*$ is the $n$th Morava $K$--theory.

To define the second family, recall that a $p$--local finite complex
$M$ is of {\em type $n$} if $M$ is $K(m)_*$--acyclic for $m<n$, but is
not $K(n)_*$--acyclic.  If $M$ is of type $n$, then $M$ admits a {\em
$v_n$--self map}, a $K(n)_*$--isomorphism
$$ v \co \Sigma^d M \ra M.$$
We let $T(n)$ denote the telescope of $v$.  A consequence of the
Nilpotence and Periodicity Theorems of Devanitz, Hopkins, and 
Smith
\cite{dhs,hs,ravenel} is that the associated Bousfield localization
functor $L_{T(n)} \co \Sp \ra \Sp$ is independent of the choice of both
the complex and self map.  Also, we recall that $T(n)_*$--acyclics are
$K(n)_*$--acyclic; thus the associated localization functors are related
by $L_{K(n)} \simeq L_{K(n)}L_{T(n)}$.

The main theorem of \cite{kuhntate} says that Goodwillie towers of
functors from spectra to spectra always split after applying $L_{T(n)}$.

\begin{thm}{\rm\cite[Theorem 1.1]{kuhntate}}\qua
\label{splitting theorem}  Let $F \co
\Sp \ra \Sp$ be any homotopy functor.  For all primes $p$, $n \geq 1$,
and $d \geq 1$, the natural map
$$ p_d(X) \co P_dF(X) \ra P_{d-1}F(X)$$
admits a natural homotopy section after applying $L_{T(n)}$. 
\end{thm}

\begin{cor} \label{tate corollary} Let $F \co \Sp \ra \Sp$ be polynomial
of degree less than $d$ and $G \co \Sp \ra \Sp$ homogeneous of degree $d$.
Then any natural transformation $f \co F \ra L_{T(n)}G$ will be null.
\end{cor}

The corollary follows formally from the theorem using \fullref{homology
error lemma}: we leave verifying this as an exercise for the reader.
The theorem is proved by combining \fullref{pullback prop} and
\fullref{homology iso prop2} with the following vanishing theorem about
Tate homology.

\begin{thm}{\rm\cite[Theorem 1.5]{kuhntate}}
\label{Tate theorem}  For all finite
groups $G$, primes $p$, and $n \geq 1$,
$$L_{T(n)}{\mathcal T}_G (L_{T(n)}S) \simeq *.$$
\end{thm}

In \cite{kuhntate}, I manage to first reduce the proof of the theorem
to the case when $G = \Sigma_p$.  There are familiar `inverse limits of
Thom spectra' models for $L_{T(n)}{\mathcal T}_{\Sigma_p} (L_{T(n)}S)$.
Using these, the equivalence $L_{T(n)}{\mathcal T}_{\Sigma_p} (L_{T(n)}S)
\simeq *$
can be shown to be equivalent to the case when $X=S$ of the following
statement about the Goodwillie tower of $\Sinfty \Oinfty X$.

\begin{thm}{\rm\cite[Theorem 3.7]{kuhntate}}\qua
\label{vanishing thm}
$$\holim_k \Sigma^k L_{T(n)}P_p^{S^{\infty}}(\Sigma^{-k}X) \ra L_{T(n)}X$$
admits a homotopy section.
\end{thm}

This theorem follows immediately from the existence of the natural section
$\eta_n(X)$ of $L_{T(n)}e_1(X) \co L_{T(n)} \Sinfty \Oinfty X \ra L_{T(n)}X$
to be discussed in the next subsection.

\begin{rem} A weaker version of \fullref{Tate theorem} with $K(n)$
replacing $T(n)$ appears in work by Greenlees, Hovey, and Sadofsky
\cite{gs,hsadofsky}, and certainly inspired my thinking, if not my proof.
\fullref{Tate theorem} when $G = \mathbb Z/2$ is equivalent to the main
theorem of Mahowald--Shick \cite{ms}.
\end{rem}

\subsection{The periodic homology of infinite loopspaces}

Using the full strength of the Periodicity Theorem, Bousfield and I have
constructed `telescopic functors' as in the next theorem.

\begin{thm} \cite{bousfield3,k2,bousfield7} \label{tel functor theorem}
For all $p$ and $n>0$, there exists a functor $\Phi_n \co \T \ra \Sp$ such
that $\Phi_n(\Oinfty X) \simeq L_{T(n)}X$.  Furthermore $\Phi_n(Z)$
is always $T(n)$--local.
\end{thm}

Some further nice properties of $\Phi_n$ will be discussed in the next
subsection: see \fullref{phi prop}.  Here we note the following corollary.

\begin{cor} After applying $L_{T(n)}$, the natural transformation
$$e_1(X) \co \Sinfty \Oinfty X \ra X$$ admits a section
$$\eta_n(X) \co L_{T(n)}X \ra L_{T(n)} \Sinfty \Oinfty X.$$
\end{cor}

The section is defined by applying $\Phi_n$ to the natural map
$$\eta(\Oinfty X) \co \Oinfty X \ra Q \Oinfty X.$$

\begin{rem} $\eta_n$ is unique up to `tower phantom' behavior in the
following sense: for all $d$, the composite
$$ L_{T(n)}X \xra{\eta_n(X)} L_{T(n)} \Sinfty \Oinfty X
\xra{L_{T(n)}e_d(X)} L_{T(n)}P_d^{S^{\infty}}(X)$$
is the unique natural section of $L_{T(n)}P_d^{S^{\infty}}(X)
\xra{L_{T(n)}p_d(X)} L_{T(n)}X$.  Here uniqueness is a consequence of
\fullref{tate corollary}.
\end{rem}

In \cite{k6}, I use $\eta_n$ to prove a splitting result in a manner
similar to \fullref{splitting example}.  The natural transformation
$$ \eta_n(X) \co X \ra L_{T(n)} \Sinfty \Oinfty X$$
induces a map of commutative augmented $L_{T(n)}S$--algebras
$$ s_n(X) \co L_{T(n)} \mathbb P(X) \ra L_{T(n)} \Sinfty (\Oinfty X)_+.$$
As in \fullref{splitting example}, $s_n(X)$ has been constructed so
that it is easy to see that $$TAQ(s_n(X)) \co L_{T(n)}X \ra L_{T(n)}X$$
is homotopic to the identity, and one learns that $s_n(X)$ induces an
equivalence of localized Goodwillie towers.  Because the towers have been
localized with respect to a nonconnected homology theory, the convergence
of localized towers is problematic.  However, one can easily deduce the
first statement of the next theorem, and starting from this, I was able
to establish the rest.

\begin{thm} {\rm \cite{k6}} \qua For all $X \in \Sp$,
$$ s_n(X)_* \co K(n)_*(\mathbb P(X)) \ra K(n)_*(\Oinfty X)$$
is monic, and fits into a chain complex of commutative $K(n)_*$--Hopf
algebras
$$ K(n)_*(\mathbb P(X)) \xra{s_{n*}} K(n)_*(\Oinfty X) \ra
\bigotimes_{j=0}^{n+1}K(n)_*(K(\pi_j(X),j)).$$
This sequence of Hopf algebras is exact if $X$ is $T(m)_*$--acyclic for
all $0 < m < n$, and only if $X$ is $K(m)_*$--acyclic for all $0 < m < n$.
\end{thm}

Note that the two acyclicity conditions on $X$ are empty if $n=1$. They
agree if $n=2$, by the truth of the Telescope Conjecture when $n=1$.

Recall that the  Telescope Conjecture asserts that a $K(n)_*$--acyclic
spectrum will always be $T(n)_*$--acyclic, and is believed to be not
true for $n \geq 2$.
Peter May has remarked that maybe the theorem could be used to disprove
it.  I end this subsection by describing one way this might go.

If $Z$ is a connected space, let $j_n(Z) \co L_{T(n)}\mathbb P(Z) \ra
L_{T(n)}\mathbb P(Z)$ be the composite
$$ L_{T(n)}\mathbb P(Z) \xra{s_n(\Sinfty Z)} L_{T(n)} \Sinfty QZ
\xra{L_{T(n)}s(Z)^{-1}}  L_{T(n)}\mathbb P(Z).$$
Here we have written $\mathbb P(Z)$ for $\mathbb P(\Sinfty Z)$.

The theorem says that if $Z$ is $T(m)_*$--acyclic for all $0 < m < n$,
and only if $Z$ is $K(m)_*$--acyclic for all $0 < m < n$, there is a
short exact sequence of $K(n)_*$--Hopf algebras
$$ K(n)_*(\mathbb P(Z)) \xra{j_{n*}} K(n)_*(\mathbb P(Z)) \ra
\bigotimes_{j=0}^{n+1}K(n)_*(K(\pi_j^S(Z),j)).$$

It appears that for some $Z$, a calculation of both $K(n)_*(\mathbb P(Z))$
and $j_{n*}$ may be accessible.  If one could find a $K(n-1)_*$--acyclic
space\footnote{If a space $Z$ is $K(n-1)_*$--acyclic, then it is
$K(m)_*$--acyclic for all $m<n$ by Bousfield \cite{bousfield99}.} $Z$, and
explicit calculation showed that the above sequence is {\em not} exact,
it would follow that $Z$ would not be $T(m)_*$--acyclic for some $0<m<n$.
The first example to check is when $p=2$, $n=3$, and $Z = K(\mathbb
Z/2,3)$: it is known that this space is $K(2)_*$--acyclic, but it is
unknown whether or not it is $T(2)_*$--acyclic.

\subsection{The periodic homotopy groups of odd dimensional spheres}

Let $v \co \Sigma^d M \ra M$ be a $K(n)_*$--isomorphism of a space $M$ whose
suspension spectrum is a finite compex of type $n$.  If $Z$ is a space,
one can use $v$ to define periodic homotopy groups by letting
$$\displaystyle v^{-1}\pi_*(Z;M) = \colim_r [\Sigma^{rd}M,Z]_*.$$
It is clear that these behave well with respect to fibration sequences
in the $Z$ variable.  These can be similarly defined for spectra, and
it is evident that there is an isomorphism $v^{-1}\pi_*(\Oinfty X;M)=
v^{-1}\pi_*(X;M)$.

The direct limit appearing in the definition suggests that these functors
of spaces do not necessarily commute with holimits of towers.  However
Arone and Mahowald note that the properties listed in \fullref{Arone
theorem} imply that the tower for an odd dimensional sphere leads to a
convergent spectral sequence with only a finite number of infinite loop
fibers for computing periodic homotopy.  More precisely, they show prove
the following.

\begin{thm} \label{AMthm} Let $m$ be odd. With $(M,v)$ as above, the
natural map
$$ v^{-1}\pi_*(S^m;M) \ra v^{-1}\pi_*(P_{p^n}S^m;M)$$
is an isomorphism.
\end{thm}

Bousfield notes that periodic homotopy can be computed using the
telescopic functor $\Phi_n$.

\begin{prop}[Bousfield {{\cite[Theorem 5.3(ii) and Corollary
5.10(ii)]{bousfield7}}}]
\label{phi prop} There are natural isomorphisms 
$$v^{-1}\pi_*(Z;M) \simeq [M,\Phi_n(Z)]_*.$$  
Furthermore, given $f \co Y \ra Z$, $v^{-1}\pi_*(f;M)$
is an isomorphism if and only if $\Phi_n(f) \co \Phi_n(Y) \ra \Phi_n(Z)$
is a weak equivalence.
\end{prop}

One can now deduce the following theorem.

\begin{thm} \label{periodic theorem 1} When $m$ is odd, there is a
spectral sequence which converges to $v^{-1}\pi_*(S^m;M)$ with
\begin{equation*}
E^1_{-k,*} =
\begin{cases}
[M,L_{T(n)}L(k,m)]_{*+k-m} & \text{for } 0 \leq k \leq n \\ 0 &
\text{otherwise.}
\end{cases}
\end{equation*}
\end{thm}

This follows by assembling various results from above.
By \fullref{AMthm} and \fullref{Arone theorem}(1), applying
$v^{-1}\pi_*(\qua;M)$ to the tower $P_d(S^m)$, the
tower of the identity evaluated at $S^m$, yields a spectral sequence
with $$E^1_{-k,*} = v^{-1}\pi_*(D_{p^k}(S^m);M)$$ for $0\leq k\leq n$
and 0 otherwise.
But this rewrites:
\begin{align*}
v^{-1}\pi_*(D_{p^k}(S^m);M) &
= [M,\Phi_n(\Omega^{\infty}D_{p^k}^{st}(S^m)]_* & \text{(by \fullref{phi
prop})}\\
  & = [M,L_{T(n)}D_{p^k}^{st}(S^m)]_*  &\text{(by \fullref{tel functor
  theorem})} \\
  & = [M,L_{T(n)}L(k,m)]_{*+k-m} &\text{(by \fullref{Arone
  theorem}(2))}.
\end{align*}
Similarly, applying $\Phi_n$ to the tower $P_d(S^m)$, yields the next
theorem.

\begin{thm} \label{phi theorem}  Let $m$ be odd.  The spectrum
$\Phi_n(S^m)$ admits a finite decreasing filtration with fibers
$L_{T(n)}\Sigma^{m-k}L(k,m)$ for $k = 0, \dots, n$.
\end{thm}

For a space $Z$, the $v_n$--periodic homotopy groups with integral
coefficients are defined as follows.

Let $L_{n-1}^f$ denote localization with respect to $T(0) \vee \dots
\vee T(n-1)$.  Then $L_{n-1}^fX \simeq *$ if and only if $X$ is the
filtered homotopy colimit of finite spectra of type at least $n$, and
it follows quite formally from this that $L_{n-1}^f$ is smashing. (See
eg Bousfield \cite{bousfield7} and Miller \cite{miller}.)

Let $C_{n-1}^fX$ denote the fiber of $X \ra L_{n-1}^fX$.  Again
$C_{n-1}^f$ is smashing. $C_{n-1}^fS$ can be realized as the homotopy
colimit of a sequence of finite spectra of type $n$,
$$ C_1 \ra C_2 \ra \cdots .$$
Fairly explicit ways of doing this are described in
Mahowald--Sadofsky \cite{mahowaldsadofsky} and \cite{k2}.

One then defines $ v_n^{-1}\pi_*(Z)$ by letting
$$ v_n^{-1}\pi_*(Z) = \colim_r v^{-1}\pi_*(Z;\D(C_r)).$$
The results above show that alternatively this can be computed as
$$ v_n^{-1}\pi_*(Z) = \pi_*(C_{n-1}^f\Phi_n(Z)).$$
A little diagram chasing shows that there is an equivalence
$$ C_{n-1}^fL_{T(n)}X \simeq M_n^fX,$$
where $M_n^fX$ is the fiber of $L_n^fX \ra L_{n-1}^fX$. ($M_n^fX$ is the
$n^{th}$ monocular part of $X$ is the terminology of 
Bousfield \cite{bousfield3}.)

Analogous to Theorems \ref{periodic theorem 1} and \ref{phi theorem}
is the following.

\begin{thm} \label{periodic theorem 2} When $m$ is odd, there is a
spectral sequence for computing $v_n^{-1}\pi_*(S^m)$ with
\begin{equation*}
E^1_{-k,*} =
\begin{cases}
\pi_{*+k-m}(M_n^fL(k,m)) & \text{for } 0 \leq k \leq n \\ 0 &
\text{otherwise.}
\end{cases}
\end{equation*}
\end{thm}

\fullref{periodic theorem 1}, \fullref{phi theorem}, and 
\fullref{periodic theorem 2} can be improved. 
We describe how this goes for \fullref{phi
theorem}.  The inclusion
$$ S^1 \ra \Omega^{m-1}S^m$$
induces a map of towers
$$ P_d(S^1) \ra \Omega^{m-1}P_d(S^m).$$
Taking homotopy fibers and localizing at $p$, one gets a tower
converging to $\Omega_0^mS^m$, with fibers $\Oinfty L(k)_1^{m-1}$,
where $L(k)_1^{m-1}$ is the fiber of the natural map of spectra $L(k,1)
\ra L(k,m)$.  Applying $\Phi_n$ to this tower, one deduces

\begin{thm} \label{phi theorem 2}  Let $m$ be odd.  The spectrum
$\Phi_n(S^m)$ admits a finite decreasing filtration with fibers
$L_{T(n)}\Sigma^{m+1-k}L(k)_1^{m-1}$ for $k = 1, \dots, n$.
\end{thm}

\begin{exm} When $p=2$, $L(1)_1^m = \mathbb RP^{m}$.  Specializing to
$n=1$,  we learn that there is a weak equivalence
$$ \Phi_1(S^{2k+1}) \simeq \Sigma^{2k+1}L_{K(1)}\mathbb RP^{2k}.$$
Specializing to $n=2$, we learn that there is a fibration sequence
of spectra
$$ \Phi_2(S^{2k+1}) \ra \Sigma^{2k+1}L_{T(2)}\mathbb RP^{2k} \ra
\Sigma^{2k}L_{T(2)} L(2)_1^{2k}.$$
The first of these is equivalent to Mahowald's theorem \cite{mahowald}
that the James--Hopf map $\Omega^{2k} S^{2k+1} \ra \Oinfty \Sinfty
\Sigma \mathbb RP^{2k}$ induces an isomorphism on $v_1$--periodic
homotopy groups.  Thus Mahowald's older work can be given a conceptual
proof and put into a general context.
\end{exm}

\bibliographystyle{gtart}
\bibliography{link}

\end{document}